\documentclass[12pt]{amsart}

\usepackage{a4}
\usepackage{amsmath}
\usepackage{amssymb}
\usepackage{amsthm}
\numberwithin{equation}{section}

\newcommand{\A }{\mathcal{A}}
\newcommand{\cA }{\mathcal{A}}

\newcommand{\B }{\mathcal{B}}

\newcommand{\cU }{\mathcal{U}}

\newcommand{\cop }{\mathrm{cop}}
\newcommand{\C }{\mathbb{C}}
\newcommand{\dif }{\mathrm{d}}

\newcommand{\eins }{{\bf 1}}
\newcommand{\End}{\mathrm{End}}

\newcommand{\id}{\textrm{id}}
\newcommand{\ids}{\mbox{\scriptsize{Id}}}
\newcommand{\im}{\textrm{Im}}

\newcommand{\Hom}{\mathrm{Hom}}
\newcommand{\kopr }{\varDelta }
\newcommand{\kow }{\varDelta }

\newcommand{\lact }{\triangleright}
\newcommand{\lid }{\mathcal{L}}
\newcommand{\Lin }{\mathrm{Lin}}
\newcommand{\N}{\mathbb{N}}
\newcommand{\ot }{\otimes }
\newcommand{\op }{\mathrm{op}}

\newcommand{\pair }[2]{\langle #1,#2\rangle }

\newcommand{\podl }{\mathcal{O}_q(\mathbb{S}^2_c)}
\newcommand{\podlo }{\mathcal{O}_q(\mathbb{S}^2_0)}

\newcommand{\rform }{\mathbf{r}}
\newcommand{\SO }[1]{\mathcal{O}_{q^2}(\mathrm{SO}(#1))}
\newcommand{\SLZ }{\mathcal{O}_q(\mathrm{SL}(2))}

\newcommand{\U }{U}

\newcommand{\Uqbp }{U_q(\mathfrak{b}_+)}
\newcommand{\Uqbm }{U_q(\mathfrak{b}_-)}
\newcommand{\Uqnp }{U_q(\mathfrak{n}_+)}
\newcommand{\Uqnm }{U_q(\mathfrak{n}_-)}
\newcommand{\UslZ }{U_q(\mathfrak{sl}_2)}

\newcommand{\vep }{\varepsilon }

\newtheorem{lemma}{Lemma}[section]
\newtheorem{cor}{Corollary}[section]
\newtheorem{theorem}{Theorem}[section]
\newtheorem{proposition}{Proposition}[section]

\theoremstyle{remark}
\newtheorem*{remark}{Remark}

\begin{document}

\title[Podle\'s' Quantum Sphere]{Podle\'s' Quantum Sphere:
  Dual Coalgebra and Classification of Covariant First Order Differential
  Calculus}

\author{Istv\'an Heckenberger}
\address{Mathematisches Institut, Universit\"at Leipzig, Augustusplatz 10,
         04109 Leipzig, Germany}

\email{heckenbe@mathematik.uni-leipzig.de\\
       kolb@itp.uni-leipzig.de} 

\author{Stefan Kolb}

\subjclass[2000]{58B32, 81R50}



\keywords{Quantum groups, Podle\'s' quantum sphere, dual coalgebra, 
  differential calculus} 

\begin{abstract}
  The dual coalgebra of Podle\'s' quantum sphere $\podl$ is determined
  explicitly. This result is used to classify all finite dimensional
  covariant first order differential calculi over $\podl$ for all but
  exceptional values of the parameter $c$.
\end{abstract}

\maketitle
\section{Introduction}
Podle\'s' quantum sphere $\podl$ \cite{a-Po87} is one of the best investigated
examples of a quantum space, i.~e.~of a comodule algebra over the
$q$-deformed coordinate ring of some affine algebraic group.
Nevertheless, classification of covariant first order differential calculus
(FODC) over $\podl$ in the sense of Woronowicz \cite{a-Woro2} has so far
been achieved only under
additional assumptions and in low dimensions. In \cite{a-Po92} certain
2-dimensional covariant FODC over $\podl$ which in many respects behave
similarly as their classical counterparts have been classified. It turned out
that only in the so called quantum subgroup case $c=0$ such a calculus exists
and is then uniquely determined. 
All covariant FODC which as right modules are freely generated by the
differentials of the generators $e_i$, $i=-1,0,1$ of $\podl$ have been
determined in \cite{a-ApelSchm94}. It was shown by computer calculations
that for all but exceptional values of $c$ exactly one such calculus exists.
Finally in \cite{a-Herm01} a general notion of dimension of covariant FODC was
introduced and all 2-dimensional covariant FODC over $\podl$ have been
classified.

In the present paper \textit{all} finite dimensional covariant FODC over
$\podl$ for all but exceptional values of $c$ are classified.
It turns out, that for generic $c$ there exists precisely one
irreducible covariant FODC for any irreducible $\SLZ$-subcomodule of
$\podl$. The subcomodule $\C\cdot \eins$ corresponds to the trivial calculus
while in general the irreducible differential calculus has the same dimension
as the corresponding $\SLZ$-subcomodule. For generic $c$ any covariant FODC
over $\podl$ can be uniquely written as a direct sum of irreducible FODC.
The exceptional cases include the quantum subgroup case $c=0$.

The main tool on this way is the notion of quantum tangent space introduced
for quantum groups in \cite{a-Woro2} and generalized to a large class of
quantum spaces in \cite{a-HeKo01p}. Podle\'s' quantum sphere can be obtained
as right $K_c$-invariant elements in $\SLZ$ where $K_c$ denotes a left coideal
subalgebra of $\UslZ$ generated by one twisted primitive element $X_c$.
The notion of quantum tangent space allows one to identify finite dimensional
covariant FODC over $\podl$ with finite dimensional left subcomodules
$T_\vep\subset \podl^\circ$ of the dual coalgebra which are right
$K_c$-invariant and contain the counit $\vep$. Thus as a first step towards
classification the dual coalgebra $\podl^\circ$ is determined explicitly in
Theorem \ref{dualcoalg}. It turns out that for all but exceptional values of
$c$ the restriction $\SLZ^\circ\rightarrow \podl^\circ$ is onto.

Next, the subspace $F(\podl^\circ,K_c)$ of elements of $\podl^\circ$ with
finite right $K_c$-action is determined. The action of the generator $X_c$
induces a $\UslZ$-action on $F(\podl^\circ,K_c)$ such that the decomposition
into irreducible $\UslZ$-modules corresponds to the decomposition into right
$K_c$-invariant left $\podl^\circ$-comodules.
To calculate $F(\podl^\circ,K_c)$ explicit results of \cite{a-MullSch} are
employed.

The quantum tangent spaces of the covariant FODC constructed
in \cite{a-Herm98} are calculated. It turns out that for generic $c$ the
resulting tangent spaces cover all tangent spaces obtained in the
classification. Therefore up to exceptional values of $c$ all covariant
FODC over $\podl$ can be constructed by this method.
Moreover it is shown in Proposition \ref{freiinner} that these FODC are free
left and right $\podl$-modules and inner calculi.

The ordering of this paper is as follows. In Section \ref{podlsphere}
the definition and some properties of $\podl$ are recalled. Section 
\ref{dualcoalgebra} serves to give a complete description of the dual
coalgebra $\podl^\circ$. The main idea on this way is to show that all
representations of $\podl$ can be written as direct sums of representations
of certain localizations of $\podl$. These localizations are seen to be
isomorphic to $\Uqbm^\op$ and the dual coalgebra of $\Uqbm^\op$ is known
\cite{b-Joseph}.
In Section \ref{podlocfin} the subspace $F(\podl^\circ,K_c)$ is determined
and decomposed into $\UslZ$-modules. The notion of covariant FODC and
quantum tangent space are recalled in the last section. Combination of the
above steps lead to the classification result in Theorem \ref{podlclasstheo}. 

If not stated otherwise all notations and conventions coincide with those
introduced in \cite{b-KS}. All through this paper $q\in\C\setminus\{0\}$ will
be assumed not to be a root of unity. For any element $a$ of a coalgebra
$\A$ with counit $\vep$ define $a^+:=a-\vep(a)$ and for any subset
$\B\subset \A$ set $\B^+:=\{b^+\,|\, b\in \A\}$.

\section{Podle\'s' Quantum Sphere}\label{podlsphere}
Let $u^i_j$, $i,j=1,2,$ denote the matrix coefficients of the vector
representation of $\UslZ$, i.e. the generators of the quantum group
$\SLZ$. In the notation of \cite{a-Po87} the matrix coefficients of the three
dimensional representation of $\UslZ$ are given by
\begin{align*}
  (\pi^i_j)_{i,j=-1,0,1}=\left(\begin{array}{ccc}
                         u^2_2 u^2_2&-(q^2+1)u^2_2u^2_1&-qu^2_1u^2_1\\
                         -q^{-1}u^1_2u^2_2&1+(q+q^{-1})u^1_2u^2_1&u^1_1u^2_1\\
                         -q^{-1}u^1_2u^1_2&(q+q^{-1})u^1_2u^1_1&u^1_1u^1_1
                         \end{array}\right).
\end{align*}  
For $c=(\vep(e_{-1})\vep(e_1):\vep(e_0)^2)\in \C P^1$ Podle\'s' quantum
sphere $\B=\podl$ \cite{a-Po87} is isomorphic to the subalgebra of $\SLZ$
generated by $e_i=\sum_{j=-1,0,1}\vep(e_j)\pi^j_i$, $i=-1,0,1$.
The algebra $\podl$ obtains the structure of a right $\SLZ$-comodule algebra
by $\kow(e_i)=e_j\ot\pi^j_i$.
A complete set of defining relations of $\podl$ is given by
\begin{align*}
  (1+q^2)(e_{-1}e_1+q^{-2}e_1e_{-1})+e_0^2&=\rho\\
  -q^2e_{-1}e_0+e_0e_{-1}&=\lambda e_{-1}\\
  (1+q^2)(e_{-1}e_1-e_1e_{-1})+(1-q^2)e_0^2&=\lambda e_0\\
  e_1e_0-q^2e_0e_1&=\lambda e_1
\end{align*}
where $\rho=q^{-2}(q^2+1)^2\vep(e_{-1})\vep(e_1)+\vep(e_0)^2$ and
$\lambda=(1-q^2)\vep(e_0)$. For $c\neq \infty$ one can choose $\vep(e_0)=1$,
$\vep(e_{-1})=\vep(e_1)$. Then $\lambda=1-q^2$ and $\rho=(q+q^{-1})^2c+1$.
Defining $A=(1+q^2)^{-1}(1-e_0)$ the above relations can be rewritten as
\begin{align}
e_{-1}e_1&=A-A^2+c\label{e-e}\\
e_1e_{-1}&=q^2A-q^4A^2+c\label{ee-}\\
e_1 A&=q^2 A e_{1}\label{eA}\\
e_{-1} A &=q^{-2} A e_{-1}\label{e-A}.
\end{align}  
Similarly for $c=\infty$ choose $\vep({e_0})=0$ and $\vep(e_{-1})=\vep(e_1)=1$,
i.e. $\lambda=0$ and $\rho=(q+q^{-1})^2$. Defining $A=-(1+q^2)^{-1}e_0$
the above relations are equivalent to
\begin{align}
e_{-1}e_1&=-A^2+1\label{e-ei}\\
e_1e_{-1}&=-q^4A^2+1\label{ee-i}\\
e_1 A&=q^2 A e_{1}\label{eAi}\\
e_{-1} A &=q^{-2} A e_{-1}\label{e-Ai}.
\end{align}
Define linear functionals $f_\lambda$, $\lambda\in \C\setminus\{0\}$ and $g$
in the dual Hopf algebra $\SLZ^\circ$ of $\SLZ$ by
\begin{align}\label{fugu}
  f_\lambda((u^i_j))&=\begin{pmatrix}\lambda&0\\0&\lambda^{-1}\end{pmatrix},&
  g((u^i_j))&=\begin{pmatrix}1&0\\0&-1\end{pmatrix}
\end{align}
and
\begin{align}\label{koprfg}
  \kopr f_\lambda&=f_\lambda\ot f_\lambda, &\kopr g&=g\ot \vep+\vep\ot g.
\end{align}
Recall that the dual pairing \cite{b-KS} between $\UslZ$ and $\SLZ$ induces
linear functionals $E$ and $F$ in $\SLZ^\circ$ satisfying
\begin{align}\label{EuFu}
E((u^i_j))&=\begin{pmatrix}0&0\\1&0\end{pmatrix},&
F((u^i_j))&=\begin{pmatrix}0&1\\0&0\end{pmatrix}
\end{align}
and
\begin{align}\label{koprEF}
\kopr E &= E \ot K+\vep\ot E,  & \kopr F &= F \ot \vep + K^{-1}\ot F
\end{align}
where $K=f_{q^{-1}}$.
Let $\cU\subset \SLZ^\circ$ denote the algebra generated by the functionals
$f_\lambda$, $\lambda\in \C\setminus\{0\}$, $E$, $F$ and $g$.
For transcendental $q$ the Hopf algebra $\cU$ is isomorphic to $\SLZ^\circ$
\cite[9.4.9]{b-Joseph}. The above functionals satisfy the relations
\begin{align}\label{Urel}
\begin{aligned}  
  f_\lambda f_\mu&=f_{\lambda\mu},& f_\lambda E&=\lambda^{-2}E f_\lambda,&
                       f_\lambda F&=\lambda^{2}F f_\lambda,\\
  f_\lambda g&=g f_\lambda,& Eg&=(g+2)E,&Fg&=(g-2)F,\\
   &&EF-FE&=\frac{K-K^{-1}}{q-q^{-1}}.&&
\end{aligned}
\end{align}

Note that the subalgebra of $\SLZ^\circ$ generated by $E,F,K$ and $K^{-1}$ is
isomorphic to $\UslZ$, \cite[4.4.1]{b-KS}.
Evaluating the functionals $f_\lambda,g,E$ and $F$ on the matrix
coefficients $\pi^i_j$ one obtains
\begin{align}
  \begin{aligned}  
    f_\lambda((\pi^i_j))&=\lambda^{-2} E^{-1}_{-1}+E_0^0+\lambda^2 E^1_1\\
    g((\pi^i_j))&=-2 E^{-1}_{-1}+ 2E^1_1\\
    E((\pi^i_j))&=-(q^2+1) E^{-1}_0+E^0_1\\
    F((\pi^i_j))&=-q^{-1} E^{0}_{-1}+(q+q^{-1})E^1_0
  \end{aligned}
\end{align}  
where $E^i_j$, $i,j=-1,0,1$ denotes the $3\times 3$-matrix with entry $1$
at position $(i,j)$ and zero elsewhere.
The right comodule structure of $\podl$ induces a left action of $\cU$ on
$\podl$ which is given by
\begin{align}
f_\lambda \lact e_i&=\lambda^{2i}e_i\label{lambei}\\
g\lact e_i&=2i e_i\label{gei}\\
E\lact e_i&=-(q^2+1)\delta_{i,0}e_{-1} +\delta_{i,1}e_{0}\label{Eei}\\
F\lact e_i&=-q^{-1}\delta_{i,-1}e_{0} +(q+q^{-1})\delta_{i,0}e_{1}\label{Fei}.
\end{align}
For $n\in\N_0/2$ set $c(n)=-1/(q^n+q^{-n})^2$.
Since $q$ is not a root of unity $c(n)\neq c(m)$ for all $n,m\in\N_0/2$,
$n\neq m$.
Define subsets of $\C P^1$ by
\begin{align*}
J_1&:=\{c\in\C P^1\,|\,c\neq c(n)\quad \forall n\in\N/2\setminus \N\}\\
J_2&:=\{c\in\C P^1\,|\,c\neq c(n)\quad \forall n\in\N_0/2\}\\
\end{align*}
It is known \cite[Rem.~4.5.3]{a-MullSch} that the following statements are
equivalent:
\begin{enumerate}
  \item $c\in J_1$ 
  \item $\podl\cong \{b\in\SLZ\,|\,X(b_{(1)}) b_{(2)}=0\}$ for a twisted
      primitive element
      \begin{align*}
         X=\alpha (K^{-1}-1)+\beta K^{-1}E + \gamma F\in\cU
      \end{align*}
      and
      \begin{equation*}
         c=\begin{cases}\frac{\beta\gamma q^{-1}}{\alpha^2(q-q^{-1})^2}
            &\textrm{if }\alpha\neq 0,\\
           \infty&\textrm{if }\alpha=0 \textrm{ and } \beta\gamma\neq 0.
           \end{cases}
      \end{equation*}
\end{enumerate}      
Calculating the pairing between $X$ and the explicit generators $e_i\in\SLZ$
chosen above one obtains $-\vep(e_1)(q-q^{-1})\alpha=\gamma$ and
$\beta=q\gamma$ in the case $c\neq\infty$. Similarly for $c=\infty$
one obtains $\alpha=0$ and $\beta=q\gamma$.
Thus the embeddings from above are realized by
\begin{align}\label{Xc}
X_c=\begin{cases}qK^{-1}E+F& \textrm{if } c=\infty,\\
                 K         & \textrm{if } c=0,\\
                 -(c^{1/2}(q-q^{-1}))^{-1}(K^{-1}-1)+qK^{-1}E+F&\textrm{else } 
    \end{cases}             
\end{align}  
for any square root $\vep({e_1})=c^{1/2}$ of $c$.
Define $K_c=\C[X_c]\subset \UslZ$.
If $c\in J_2$ then any finite dimensional
$\UslZ$-module is a direct sum of irreducible $K_c$-modules and therefore
$\SLZ$ is a faithfully flat
left (and right) $\podl$-module \cite[Thm.~5.2]{a-MullSch}.

\section{The Dual Coalgebra $\B^\circ=\podl^\circ$}\label{dualcoalgebra}
To understand the dual coalgebra \cite[Sect.~6.0]{b-Sweedler} $\podl^\circ$
of Podle\'s' quantum sphere it is useful to consider first the dual Hopf
algebra $(\Uqbm^{\textrm{op}})^\circ$ where
$\Uqbm\subset \UslZ$ denotes the subalgebra generated by $F,K$ and $K^{-1}$.
Further let $U_0$, $\Uqnp$, $\Uqnm$ and $\Uqbp$ denote the subalgebra of
$\UslZ$ generated by $\{K,K^{-1}\}$, $E$, $F$ and $\{E,K,K^{-1}\}$,
respectively. By \cite[Thm. 2.1.8]{b-Joseph} the dual Hopf algebra
$(U_0)^\circ$ is isomorphic to the commutative Hopf algebra
\begin{align*}
  \C[\gamma,\chi_\lambda\,|\,\lambda\in\C\setminus\{0\}]\big/(\chi_\lambda
  \chi_\mu=\chi_{\lambda\mu},\,\chi_1=1)
\end{align*}
where $\gamma(K)=1$, $\chi_\lambda(K)=\lambda$ and the coalgebra structure
is given by
\begin{align}\label{koalg}
 \begin{aligned}\kow \gamma&=\gamma\ot 1+1\ot \gamma \\
             \kow \chi_\lambda&=\chi_\lambda\ot \chi_\lambda.
 \end{aligned}            
\end{align}  
The subalgebra $\Uqnp\subset\UslZ$ is a right $U_0$-comodule with coaction
\begin{align*}
  \delta_R(E^i)=E^i\ot K^{-i}
\end{align*}  
and therefore has a left $(U_0)^\circ$-module structure. 
The corresponding left crossed product algebra $\Uqnp\rtimes (U_0)^\circ$ 
is a Hopf algebra with $\kow E=1\ot E+E\ot \chi_{q^{-2}}$ containing
$\Uqbp$ where $K\in\Uqbp$ corresponds to $\chi_{q^{-2}}$.
The dual pairing of Hopf algebras (in the conventions of \cite[6.3.1]{b-KS})
\begin{align}\label{pair}
  \pair{\cdot}{\cdot}:\Uqbp\ot \Uqbm^{\textrm{op}}\rightarrow \C
\end{align}
given by $\pair{K}{K}=q^{-2}$, $\pair{K}{F}=\pair{E}{K}=0$ and
$\pair{E}{F}=1/(q^{-1}{-}q)$ extends to a pairing of Hopf algebras
\begin{align}\label{extpair}
  \pair{\cdot}{\cdot}:\left(\Uqnp\rtimes (U_0)^\circ\right)\ot
  \Uqbm^{\textrm{op}}\rightarrow \C
\end{align}
such that
\begin{align*}
  \pair{\gamma}{K}=1,\quad \pair{\chi_\lambda}{K}=\lambda,\quad
  \pair{\gamma}{F}=\pair{\chi_\lambda}{F}=0.
\end{align*}  
\begin{lemma}\label{paircalc}
  For $a\in \Uqnp$, $u\in U_0$, $b\in \Uqnm$ and $f\in (U_0)^\circ$ one has
  \begin{align*}
    \langle a f,b u\rangle =f(u)\langle a,b\rangle.
  \end{align*}  
  In particular the pairing (\ref{extpair}) is non-degenerate.
\end{lemma}  
\begin{proof}
Note first that $\pair{f}{b} =\vep(b)f(1)$ and therefore
\begin{align*}
\pair{f}{bu}=\pair{f_{(1)}}{u}\pair{f_{(2)}}{b}=\vep(b)\pair{f}{u}.
\end{align*}
Using this relation one calculates
\begin{align*}
 \pair{af}{bu}=\pair{a}{b_{(1)}u_{(1)}} \pair{f}{b_{(2)}u_{(2)}}
              =\pair{a}{bu_{(1)}}\pair{f}{u_{(2)}}=\pair{a}{b}\pair{f}{u}
\end{align*}
where in the last equation the property
$\pair{a}{bu}=\vep(u)\pair{a}{b}$ of (\ref{pair}) is used.
The non-degeneracy of (\ref{extpair}) now follows from the non-degeneracy of
(\ref{pair}).
\end{proof}  
By the above lemma the map of Hopf algebras
\begin{align}
  \Phi:\left(\Uqnp\rtimes (U_0)^\circ\right)\rightarrow
  (\Uqbm^{\textrm{op}})^\circ
\end{align}  
induced by (\ref{extpair}) is injective.
The following result is proven in {\em\cite[9.4.8]{b-Joseph}} for
transcendental $q$. Yet it also holds for $q\in\C\setminus\{0\}$ not a root
of unity and is reproduced here in our setting for the convenience of the
reader.
\begin{proposition} \label{uqbmdual}
  The map $\Phi$ is an isomorphism.
\end{proposition}  
\begin{proof}
  Recall that there is a canonical isomorphism $\Uqbm^\op\cong \Uqnm\ot U_0$
  of vector spaces.
  Let $J\subset \Uqbm^{\textrm{op}}$ denote any two sided ideal of finite
  codimension.
  Then $J$ contains some ideal
  $I\subset U_0$ of finite codimension and $(\Uqnm^+)^n$ for some $n\in \N$.
  Therefore $J$ contains the left ideal
  \begin{align*}
    (\Uqnm^+)^n\ot U_0 +\Uqnm\ot I \subset \Uqnm\ot U_0
  \end{align*}
  of finite codimension. Thus
  \begin{align*}
    \left(\Uqbm^{\textrm{op}}\big/ J\right)^\circ
    &\subset \left((\Uqnm\big/(\Uqnm^+)^n) \ot (U_0/I)\right)^*\\
    &=\left(\Uqnm\big/(\Uqnm^+)^n\right)^* \ot (U_0/I)^*\\
    & \subset \Uqnp\ot (U_0)^\circ
  \end{align*}
  where in the last inclusion one uses that $\Uqnp$ is the graded dual of
  $\Uqnm$ via the pairing (\ref{pair}). By Lemma \ref{paircalc} one obtains
  $\Uqnp\ot (U_0)^\circ\subset \im \,\Phi$ and therefore $\Phi$ is onto.
\end{proof}  

For the computation of $\podl^\circ$ some results about the representation
theory of the algebra $\podl$ are collected.

\begin{lemma}\label{rep}
  Any finite dimensional representation $\mu:\podl\rightarrow \End(V)$
  is a direct sum $\mu=\mu_0\oplus \mu_{\neq 0}$ where $\mu_0(A)$ is
  nilpotent and $\mu_{\neq 0}(A)$ is invertible.
  In particular the coalgebra $\podl^\circ$ is a direct sum $C_0\oplus
  C_{\neq 0}$ where $C_0$ and $C_{\neq 0}$ denote the coalgebras of
  matrix coefficients of finite dimensional representations of $\podl$ with
  nilpotent and invertible $A$ action, respectively.
  In addition
  \begin{enumerate}
    \item if $c\neq c(n)$ for all $n\in \N$ then $\mu_{\neq 0}=0$.
    \item if $c\neq 0$ then $\mu_0(e_{\pm 1})$ are isomorphisms.
    \item if $c=c(n)$ for some $n\in \N$ then there exists exactly one
          indecomposable representation $\mu_n:\podl\rightarrow \End(V)$
          such that $\mu_n(A)$ is invertible.
          This representation is $n$-dimensional.
    \item if $c=0$ then $C_0=C_{0+}\oplus C_{00}\oplus C_{0-}$ where
          $C_{0\pm}$ (resp. $C_{00}$) denotes the coalgebra of matrix
          coefficients of finite dimensional representations with
          invertible action of $e_{\pm 1}$ (with nilpotent action of $e_1$
          and $e_{-1}$).
  \end{enumerate}
\end{lemma}  

\begin{proof}
Relation (\ref{eA}) and (\ref{e-A}) imply that $e_1$ and $e_{-1}$ transform
the generalized eigenspace $V_\lambda$ of $A$ with corresponding eigenvalue
$\lambda$ to the generalized eigenspace $V_{q^{-2}\lambda}$ and
$V_{q^2\lambda}$, respectively.
Set $V_{\neq 0}:=\oplus_{\lambda\neq 0}V_\lambda$. Then
$V=V_0\oplus V_{\neq 0}$ is a direct sum of representations of $\podl$.

Since $q$ is not a root of unity $e_1$ and $e_{-1}$ act nilpotently on
$V_{\neq 0}$. Assume that $v\in V_{\neq 0}$ is an eigenvector of
$A$ with eigenvalue $\lambda$ such that $e_{-1}v=0$, $e_1^n v=0$ and
$w:=e_1^{n-1}v\neq 0$. Then relations (\ref{ee-}), (\ref{ee-i}) and
(\ref{e-e}), (\ref{e-ei}) applied
to $v$ and $w$ respectively imply
\begin{align}
&\begin{aligned}\label{c=cn}
  0&=q^2\lambda-q^4\lambda^2+c\\
  0&=q^{-2(n-1)}\lambda-q^{-4(n-1)}\lambda^2+c
\end{aligned}&\textrm{ for } c\neq \infty\\[.5\baselineskip]
&\begin{aligned}
  0&=-q^4\lambda^2+1\\
  0&=-q^{-4(n-1)}\lambda^2+1
\end{aligned}&\textrm { for } c=\infty.
\end{align}
The second set of equations cannot be fulfilled as $q$ is not a root of unity.
The first set of equations implies $c=c(n)$ and therefore proves 1).

Since $\mu_0(A)$ is nilpotent the second statement follows from
(\ref{e-e}) and (\ref{e-ei}).

To prove the third statement assume first that there exists $u\in V_{\neq 0}$
such that $(A{-}\nu)^2 u=0$ but $(A{-}\nu)u\neq 0$ for some
$\nu \in \C\setminus\{0\}$. Applying $e_{-1}$ several
times we may assume using the notations from above that $\nu=\lambda$ and
$(A{-}\nu)u=v$. Then (\ref{c=cn}) implies $\lambda=q^{n-2}/(q^n+q^{-n})$.
The relation $e_{-1}v=0$ implies that $e_{-1}u$ is an eigenvector of $A$
with corresponding eigenvalue $q^2\lambda$ or $e_{-1}u=0$. Suppose that
$e_{-1}^ku=0$ for some $k\ge 1$ and $e_{-1}^{k-1}u\neq 0$. Then on the one
hand the eigenvalue of $A$ corresponding to $e_{-1}^{k-1}u$ coincides with
$q^{2k-2}\lambda$ while on the other hand by (\ref{c=cn}) it is equal to
$\lambda$. Therefore $k=1$ and $e_{-1}u=0$.
By equation (\ref{ee-}) and $(A-\lambda)^2u=0$ one now obtains
\begin{equation*}
  (q^2-2q^4\lambda)Au+(c+q^4\lambda^2)u=-q^2\frac{q^n-q^{-n}}{q^n+q^{-n}}Au+
  (c+q^4\lambda^2)u=0.
\end{equation*}  
As $n\ge 1$ and $q^{2n}\neq 1$ this is a contradiction to the assumption that
$u$ is not an eigenvector of $A$. Thus $A$ is diagonalisable.
The relations (\ref{c=cn}) imply that all eigenvalues of $A$ lie in the set
$\{q^{n-2k}/(q^n+q^{-n})\,|\,k=1,2,\dots,n\}$. In view of (\ref{e-e})
and (\ref{ee-}) the eigenspaces for different eigenvalues are isomorphic
and $V_{\neq 0}$ is the direct sum of the $\podl$-orbits of the elements of an
arbitrary basis of $V_\lambda$. These orbits have dimension $n$.

To validate the last statement note first that any finite dimensional
representation
$\mu:\podlo\rightarrow \End(V)$ is a direct sum $\mu=\mu_+\oplus\mu'$,
$V=V_+\oplus V'$ where $\mu_+(e_1)$ is invertible and $\mu'(e_1)$ is
nilpotent. Indeed, (\ref{eA}) implies that $AV_+\subset V_+$ and
$AV'\subset V'$. On the other hand (\ref{e-e}) leads to
\begin{align*}
  e_{-1} V_+=e_{-1}e_1 V_+=(A-A^2)V_+\subset V_+
\end{align*}  
and $e_1^kV'=0$ yields
\begin{align*}
  e_1^{k+1}e_{-1}V'=e_1^k(q^2A-q^4A^2)V'\subset e_1^kV'=0.
\end{align*}  
Note then that (\ref{e-e}), (\ref{eA}) and the nilpotency of $\mu(A)$ imply
that $\mu_+(e_{-1})$ is nilpotent. Similarly $\mu'=\mu'_0\oplus\mu_-$ where
$\mu'_0(e_{-1})$ is nilpotent and $\mu_-(e_{-1})$ is invertible. 
\end{proof}  

The inclusion $\podl\subset \SLZ$ of right $\SLZ$-comodule algebras
induces a map of right $\SLZ^\circ$-module
coalgebras $\SLZ^\circ\rightarrow\podl^\circ$.
For $m,l\in \N_0$ and $\lambda \in \C\setminus\{0\}$ let
$\psi^{ml}_{\lambda^2}$
denote the image of $f_\lambda g^mE^l$ under this projection. It follows from
(\ref{lambei}) that $f_{\lambda}=f_{-\lambda}$ on $\podl$ and therefore
the definition of $\psi^{ml}_{\mu}$ does not depend on the choice of a root
of $\mu$.

\begin{theorem}\label{dualcoalg}
  The following sets form a vector space basis of $\podl^\circ$.
  \begin{enumerate}
    \item If $c\notin\{0,c(n)\,|\, n\in\N\}${\upshape:}
      $\{\psi^{ml}_\lambda\,|\,\lambda\in \C\setminus\{0\},m,l\in \N_0\}$.
    \item If $c=c(n)$, $n\in\N${\upshape:}
      $\{\psi^{ml}_\lambda\,|\,\lambda\in \C\setminus\{0\},m,l\in \N_0\}\cup
       \mathfrak{B}_n$, where $\mathfrak{B}_n$ denotes any basis of  the
       $n^2$-dimensional subspace $C_{\neq 0}$ of $\podl^\circ$.
    \item If $c=0${\upshape:} $\{E^kF^l\,|\, k,l\in \N_0\}\cup
     \{\chi_\lambda^+g^mF^l,\chi_\lambda^-g^mE^l\,|\,
       \lambda\in \C\setminus\{0\},l,m\in \N_0\}$
     where $\chi^\pm_\lambda$ is the character on $\podlo$ defined by
     $\chi_\lambda^\pm(e_i)=\delta_{i0}+\delta_{i,\pm1}\lambda^{\pm 1}$.
  \end{enumerate}
\end{theorem}  

\begin{proof}
  Consider the Hopf subalgebra $\SO3\subset \SLZ$ generated by the matrix
  coefficients $\{\pi^i_j\,|\,i,j=-1,0,1\}$ and let $J$ denote the
  intersection of the two-sided ideal $(u^1_2)\subset \SLZ$ with
  $\SO3$. There is an isomorphism of Hopf algebras
  $\SO3/J\rightarrow \Uqbm^\op$
  \begin{align*}
    u^2_2u^2_2\mapsto K^{-1},\quad u^2_1u^2_2\mapsto (1-q^2)F,\quad
    u^1_1u^1_1\mapsto K
  \end{align*}
  such that the functionals $E,f_\lambda,g\in \SO3^\circ$ given by
  (\ref{fugu}) and (\ref{EuFu}) correspond to  $E,\chi_{\lambda^2},2\gamma
  \in \Uqnp\rtimes (U_0)^\circ=(\Uqbm^\op)^\circ$.

  For $c\neq 0$ the sequence
  \begin{align*}
    \podl\hookrightarrow \SO3\rightarrow \SO3/J\to\Uqbm^\op
  \end{align*}  
  induces an isomorphism $\podl(e_{-1})\rightarrow \Uqbm^\op$
  \begin{align}\label{podbm}
  \begin{gathered}
     e_{-1}\mapsto \vep(e_{-1}) K^{-1},\quad
      e_0\mapsto \vep(e_{-1}) (q^3-q^{-1})F + \vep(e_0),\\
     e_1\mapsto -\vep(e_{-1})(q-q^{-1})^2
                 KF^2-\vep(e_0)(q-q^{-1})KF+\vep(e_1)K
  \end{gathered}
  \end{align}
  where $\podl(e_{-1})$ denotes the localization of $\podl$ with respect
  to the left and right Ore set $\{e_{-1}^n\,|\,n\in \N_0\}$.
  Thus by Lemma \ref{rep}(2) and Proposition \ref{uqbmdual} one obtains
  \begin{align*}  
     C_0\cong \podl(e_{-1})^\circ\cong (\Uqbm^\op)^\circ
          \cong \Uqnp\rtimes(U_0)^\circ
  \end{align*}
  and the basis element $\chi_\lambda\gamma^m E^l\in \Uqnp\rtimes
  (U_0)^\circ$ corresponds to $(1/2)^m\psi^{ml}_\lambda$. This proves
  (1) and one obtains (2) taking into account that for $c=c(n)$ the
  representation $\mu_{\neq 0}$ is irreducible.

  In the case $c=0$ consider the embedding
  \begin{align*}
     \podlo\hookrightarrow \SO3,\quad e_i\mapsto \pi^{-1}_i+\pi^0_i.
  \end{align*}
  Similarly to the case $c\neq 0$ this induces an isomorphism
  $\podlo(e_{-1})\rightarrow \Uqbm^\op$ given by (\ref{podbm}) with
  $\vep(e_i)=\delta_{i0}+\delta_{i,-1}$.
  Thus by Lemma \ref{rep}
  \begin{align*}
    C_{0-}\cong(\Uqbm^\op)^\circ\cong \Uqnp\rtimes (U_0)^\circ
  \end{align*}
  and  the basis element $\chi_\lambda \gamma^m E^l\in
  \Uqnp\rtimes (U_0)^\circ$ corresponds to  $(1/2)^m\chi^-_\lambda g^m E^l$.
  The subcoalgebra $C_{0+}$ is dealt with analogously replacing  the
  two-sided ideal $(u^1_2)$ by $(u^2_1)$ and $\Uqbm^\op$ by $\Uqbp^\cop$.
  The component $C_{00}$ has been shown to coincide with
  $\UslZ/(K-1)\UslZ$ in \cite[Lem.~5.2, Cor.~3.8]{a-HeKo01p}. The elements
  $\{E^kF^l\,|\,k,l\in\N_0\}$ form a basis of the coalgebra
  $\UslZ/(K-1)\UslZ$. 
\end{proof}  

\section{Local Finiteness for the $K_c$-Action on $\podl^\circ$}
    \label{podlocfin}
From now on and for the rest of this paper assume that $0\neq c\in J_2$.
For $\B=\podl$ recall that $\B^\circ$ is a right $\cU$-module.
Define  
\begin{align*}
F(\B^\circ,K_c)=\{f\in \B^\circ\,|\, \dim (fK_c)< \infty\}.
\end{align*}  
If $f\in\B^\circ$ is the restriction of an element $f'\in \cU$ to $\B$
and $k\in K_c$ then
\begin{align*}
  fk=k_{(0)}S^{-1}(k_{(-1)})f'k_{(-2)}|_\B=S^{-1}(k_{(0)})f'k_{(-1)}|_\B
\end{align*}
as $K_c$ is a left $\cU$-comodule and $k|_\B=k(1)\vep$.
Thus $F(\cU)|_{\B}\subset F(\B^\circ,K_c)$ where for any Hopf algebra
$A$
\begin{align*}
  F(A)=\{a\in A\,|\, \dim(\textrm{ad} A)a <\infty\}, \quad
  (\textrm{ad} b)a=b_{(1)}a S(b_{(2)}).
\end{align*}  

\begin{lemma}\label{fcumod}
The vector space $F(B^\circ,K_c)$ is a right $F(\cU)$-module 
left $B^\circ$-comodule. Any element of $F(B^\circ,K_c)$ is contained in a
finite dimensional right $K_c$-submodule left $B^\circ$-subcomodule. 
\end{lemma}
\begin{proof}
For $f\in F(\cU)$, $u\in F(B^\circ,K_c)$ consider $V=uK_c$ and
$W=(\textrm{ad}\,\cU)f$. Then for any $k\in K_c$
\begin{align*}
  (u\cdot f)k=uk_{(0)}\cdot S^{-1}(k_{(-1)})fk_{(-2)}\in V\cdot W.
\end{align*}
Therefore $F(B^\circ,K_c)$ is a right $F(\cU)$-module. 

Let $\bar{V}$ denote the left $\B^\circ$-comodule generated by $V$.
The vector space $\bar{V}$ is finite dimensional.
Applying the coaction to the second factor of
$k_{(-1)}\otimes u k_{(0)}\in \cU\ot V$ one obtains
\begin{align*}
  k_{-2}\ot u_{(1)}k_{(-1)}\ot u_{(2)}k_{(0)} \in\cU\ot\B^\circ\ot \bar{V}
\end{align*}  
and therefore
\begin{align*}
  u_{(1)}\ot u_{(2)}k=u_{(1)}k_{(-1)}S^{-1}(k_{(-2)})\ot u_{(2)}k_{(0)}\in
  \B^\circ\ot \bar{V}.
\end{align*}
Thus $F(\B^\circ,K_c)$ is a left $\B^\circ$-comodule and $u\in\bar{V}\supset
 \bar{V}K_c$.
\end{proof}  

\begin{lemma}\label{grad}
Any left $\B^\circ$ subcomodule $W\subset F(\B^\circ,K_c)$ is a
$(\C{\setminus}\{0\})$-graded vector space where
\begin{align*}
  \deg(g^mf_\mu E^l)=\mu.
\end{align*}  
\end{lemma}
\begin{proof}
Consider an arbitrary element $u\in W\subset F(\B^\circ,K_c)$.
By Theorem \ref{dualcoalg}.1 one can assume that
$u=\sum_{\mu} f_\mu a^\mu$ for some $a^\mu$ which are linear combinations
of basis vectors $g^mE^l$, $m,l\in \N_0$. By the explicit form (\ref{koprfg}),
(\ref{koprEF}) of the coproduct of $g$ and $E$ and by Theorem \ref{dualcoalg}.1
one can write
\begin{align*}
\kow u=\sum_{\mu} f_\mu\ot f_\mu a^\mu+ \sum_i u_i^1\ot u_i^2
\end{align*}
where $\{u_i^1,f_\mu\}$ is a set of linear independent elements in
$\B^\circ$. As $W$ is a left $\B^\circ$-comodule
$f_\mu a^\mu\in W$ for all $\mu$.
\end{proof}  
Let $F_\mu(\B^\circ,K_c)$ denote the subspace of elements of degree $\mu$ in
$F(\B^\circ,K_c)$.

\begin{lemma}
  $  F(\B^\circ,K_c)\subset \tilde{F}:=\Lin_\C \{\psi^{0l}_\lambda\,|\, l\in
  \N_0, \lambda \in\C\setminus\{0\}\}.$
  
\end{lemma}  
\begin{proof}
Consider an arbitrary element $u\in F_\mu(\B^\circ,K_c)$.
By Theorem \ref{dualcoalg}.1 one can assume that
$u=\sum_{i=0}^m g^i a_i$ for some $a_i\in \tilde{F}$ such that $\deg(a_i)=\mu$
and $a_m\neq 0$. Suppose that $m\ge 1$. Applying the coaction to $u$ one
obtains 
\begin{align*}
  \kow u= f_\mu g^{m-1}\ot (ma_m g+a_{m-1}) +\sum_i u_i^1\ot u_i^2
\end{align*}
where $\{f_\mu g^{m-1},u_i^1\}$ is a linearly independent set of elements of
$\B^\circ$. Thus, as $F(\B^\circ,K_c)$ is a left $\B^\circ$-comodule,
we can assume $m=1$. By  similar arguments one can assume that
$u=gf_\mu+a_0$ and also $f_\mu\in F(\B^\circ,K_c)$. 

One checks by direct computation that $(\textrm{ad}\,\cU)E$ is a three
dimensional vector space and therefore $E\in F(\cU)$. By Lemma \ref{fcumod}
this implies $f_\mu E^m\in F(\B^\circ,K_c)$ for all $m\in \N_0$.
Thus $g f_\mu\in F(\B^\circ,K_c)$.

Direct calculation using (\ref{Urel}) and (\ref{Xc}) leads to
\begin{align}\label{gfEXc}
g f_\mu E^l X_c&=q(q^{2l}-\mu^4)g f_{q\mu}E^{l+1}-4q\mu^4f_{q\mu}E^{l+1}+
                 \sum_{i=0}^l a_iE^i
\end{align}
where $a_i\in\Lin_\C\{gf_\nu,f_\nu\,|\,\nu\in \C\setminus\{0\}\}$.
Further
\begin{align}\label{fEXc}
\begin{aligned}  
f_\mu E^l X_c=&q(q^{2l}-\mu^4)f_{q\mu}E^{l+1}+\alpha(q^{2l}-\mu^2)f_{q\mu}E^l
               +\alpha(\mu^2-1)f_\mu E^l\\
              & + [l]\frac{q^{-l+1}K-q^{l-1}K^{-1}}{q-q^{-1}}f_\mu E^{l-1}
\end{aligned}              
\end{align}  
where as in (\ref{Xc}) $\alpha=0$ if $c=\infty$ and
$\alpha=-(c^{1/2}(q-q^{-1}))^{-1}$ else.
By (\ref{gfEXc}) 
\begin{align*}
 gf_\mu (X_c)^k=&q^k\left(\prod_{i=0}^{k-1}(q^{2i}-(q^i\mu)^4)\right)g
              f_{q^k\mu}E^{k} \\
              &-4\sum_{j=0}^{k-1}(q^j\mu)^4q^k
                 \left(\prod_{i=0 \atop i\neq j}^{k-1}
                 (q^{2i}-(q^i\mu)^4)\right)f_{q^k\mu}E^{k}+...
\end{align*}  
where $...$ denotes terms containing only smaller powers of $E$.
Therefore $gf_\mu\in F(\B^\circ)$
implies $\mu^4=q^{-2(k-1)}$ for some $k\in \N$. Then for $l\ge 0$
\begin{align*}
 gf_\mu (X_c)^{k+l}=& -4(q^{k-1}\mu)^4q^{k+l}
                 \left(\prod_{i=0 \atop i\neq k-1}^{k+l-1}
                 (q^{2i}-(q^i\mu)^4)\right)f_{q^k\mu}E^{k+l}+...
\end{align*}
again up to expressions containing only smaller powers of $E$. As
\begin{align*}
  q^{2(k+l)}-(q^{k+l}\mu)^4=q^{2(k+l)}(1-q^{2(l+1)})\neq 0
  \quad\textrm{for all } l\ge 0
\end{align*}  
the coefficient of $f_{q^k\mu}E^{k+l}$ does not vanish.
This is a contradiction to the assumption $gf_\mu\in F(\B^\circ,K_c)$. 
\end{proof}  

To shorten notation let $\psi^l_\lambda$ denote the basis element
$\psi^{0l}_\lambda$ of $\tilde{F}$. Define three maps $\phi,\varphi,\kappa:
\tilde{F}\rightarrow\tilde{F}$ by
\begin{align}
  \phi(\psi^l_\lambda)&=-\frac{q^l[l]}{q-q^{-1}}\psi^{l-1}_{q^2\lambda}
                        +\alpha q(q^{2l}-\lambda)\psi^l_{q^2\lambda}+
                        q^2(q^{2l}-\lambda^2)\psi^{l+1}_{q^2\lambda}\nonumber\\
  \varphi(\psi^l_\lambda)&=\lambda^{-1}\frac{q^{1-l}[l]}{q-q^{-1}}
                           \psi^{l-1}_{q^{-2}\lambda}\label{usl2}\\
  \kappa(\psi^l_\lambda)&=\lambda \psi^l_\lambda.\nonumber
\end{align}
In view of (\ref{fEXc}) this means
\begin{align}\label{psisummanden}
  \psi^l_\lambda X_c=q^{-1}\phi(\psi^l_\lambda)
                     +\lambda\varphi(\psi^l_\lambda)
                     +\alpha(1-\lambda^{-1})\kappa(\psi^l_\lambda).
\end{align}
Note that
\begin{align*}
  \phi\circ\varphi-\varphi\circ\phi=\frac{\kappa-\kappa^{-1}}{q-q^{-1}},\quad
  \kappa\circ\phi=q^2\phi \circ\kappa,\quad
  \kappa\circ\varphi=q^{-2}\varphi\circ\kappa,
\end{align*}  
i.e. the operators $\phi$, $\varphi$ and $\kappa$ yield a representation
$\rho:\UslZ\rightarrow \End({\tilde{F})}$, $\rho(E)=\phi$, $\rho(F)=\varphi$,
$\rho(K)=\kappa$.

\begin{lemma}\label{equiv}
For any finite dimensional subspace
$V\subset \tilde{F}$ the following statements are equivalent.
\begin{enumerate}
  \item $\kow V\subset \B^\circ\ot V$ and $VK_c\subset V$ 
  \item $V$ is a left $\UslZ$-module via $\rho$.
\end{enumerate}    
\end{lemma}  

\begin{proof}
(1)$\Rightarrow$(2) As in Lemma \ref{grad} one obtains that
$V$ is $(\C{\setminus}\{0\})$-graded. Then the assertion follows from
(\ref{psisummanden}).
To verify (2)$\Rightarrow$(1) note that
\begin{align*}
  \kow\psi^l_\lambda=\sum_{r=0}^l\begin{bmatrix}l\\r\end{bmatrix}
  q^{-r(l-r)}\psi^r_\lambda\ot\psi^{l-r}_{q^{-2r}\lambda}=
  \sum_{r=0}^l b_r \psi^r_\lambda\ot\varphi^r(\psi^{l}_\lambda)
\end{align*}
where $b_r\in\C$ depend on $r$ and $\lambda$ but not on $l$.
\end{proof}  

Lemma \ref{fcumod} implies that
$F(\B^\circ,K_c)$ is a $\rho$-invariant subspace of $\tilde{F}$.
Recall that an element $\psi\in V\setminus\{0\}$ is called a highest weight
vector of a $\UslZ$-module $V$ with highest weight $\lambda$ if
$K^{-1}\psi=\lambda\psi$ and $F\psi=0$.

\begin{proposition}\label{locfin}
There exists a decomposition of $\UslZ$-modules 
\begin{align*}
 F(\B^\circ,K_c)=\bigoplus_{\lambda\in J^c}V_{\lambda}
\end{align*}
such that
\begin{align*}
J^c&=\{q^{-l}\,|\,l\in 2\N_0\} \quad \textrm{ for }
               c\notin\left\{\infty,(q^r-q^{-r})^{-2}\,\big|\,
                  r\in \N/2\right\}\\
J^\infty&=\{\pm q^{-l}\,|\,l\in 2\N_0\}\\
J^{(q^r-q^{-r})^{-2}}&=
   \left\{q^{-l},-q^{-k}\,\big|\,l\in 2\N_0, k\in 2r+2\N_0 \right\},\quad
   r\in\N/2
\end{align*}
where the components $V_{\pm q^{-l}}$ are $(l{+}1)$-dimensional
irreducible $\UslZ$-modules of highest weight $\pm q^l$.
\end{proposition}  

\begin{proof}
Lemma \ref{fcumod} and Lemma \ref{equiv}, \textit{1.}$\Rightarrow$\textit{2.},
imply that the $\UslZ$-module $F(\B^\circ,K_c)$ can be written as a direct
sum
\begin{align*}
  F(\B^\circ,K_c)=\bigoplus_{\lambda\in J^c}V_{\lambda}
\end{align*}  
of finite dimensional irreducible $\UslZ$-modules. Here
$J^c\subset \C\setminus\{0\}$ denotes the subset of nonzero complex
numbers $\lambda$ such that $\phi$ operates nilpotently on
$\psi^0_\lambda=f_{\sqrt{\lambda}}$.
Indeed, by (\ref{usl2}) the set $\{\psi^0_\lambda\,|\,\lambda\in J^c\}$
is a basis of all highest weight vectors of $F(\B^\circ,K_c)$ with respect
to the $\UslZ$-module structure. It remains to show that $J^c$ is of the
form given in the proposition.

Note that $\phi^{l+1}(\psi^0_\lambda)=0$ and $\phi^l(\psi^0_\lambda)\neq 0$
imply $\lambda=\pm q^{-l}$. In this case the mapping
\begin{align*}
  \phi:\Lin_\C\{\psi^k_{q^{2l}\lambda}\,|\,k=0,\dots,l\}\rightarrow
  \Lin_\C\{\psi^k_{q^{2l+2}\lambda}\,|\,k=0,\dots,l\}
\end{align*}
is given by the matrix
\begin{align*}
  q \begin{pmatrix}-(\pm q^l-1)\alpha&-\hat{q}^{-1}q^0[1]&0&\cdots&0\\
        q(1-q^{2l})&-(\pm q^l-q^2)\alpha&-\hat{q}^{-1}q^1[2]&\ddots&\vdots\\
        0& q(q^2-q^{2l})&-(\pm q^l-q^4)\alpha&\ddots&0\\
        \vdots&\ddots&\ddots&\ddots&-\hat{q}^{-1}q^{l-1}[l]\\
        0&\cdots&0&q(q^{2(l-1)}-q^{2l})&   -(\pm q^l-q^{2l})\alpha   
     \end{pmatrix}
\end{align*}
with respect to the bases $\psi_\mu^k$ ($\mu=q^{2l}\lambda$ and
$\mu=q^{2l+2}\lambda$, respectively), $\hat{q}=q-q^{-1}$.
Recall that $\beta=q$ and $\gamma=1$. Using $q(1-q^{2k})=-\hat{q}q^{k+1}[k]$
the map $\phi$ can be written with respect to the bases
$\bar{\psi}^k_\mu:=(-\hat{q})^k q^{-(l-k)(l-k+1)/2}\psi^k_\mu$ as
\begin{align}\label{matrix}
  q^{l+1} \begin{pmatrix}(q^{-l}\mp 1)\alpha&[1]\gamma&0&\cdots&0\\
        q^{-l}[l]\beta&(q^{2-l}\mp 1)\alpha&[2]\gamma&\ddots&\vdots\\
        0& q^{2-l}[l-1]\beta&(q^{4-l}\mp 1)\alpha&\ddots&0\\
        \vdots&\ddots&\ddots&\ddots&[l]\gamma\\
        0&\cdots&0&q^{l-2}[1]\beta& (q^l\mp 1)\alpha   
     \end{pmatrix}.
\end{align}
In the case of minus signs in the diagonals this matrix is up to the overall
factor precisely the matrix $M_l$ describing the transpose of the left
action of $X_c$ on the $(l{+}1)$-dimensional irreducible $\UslZ$-module $V_l$
\cite[Sect.~4]{a-MullSch}.
By \cite[Prop.~4.2]{a-MullSch} the matrix $M_l$ is known to have $l+1$ not
necessarily distinct eigenvalues
\begin{align*}
  \rho_r=\frac{\alpha}{2}(q^r-q^{-r})^2+\frac{1}{2}(q^{2r}-q^{-2r})R,
 \qquad r\in I_l=\{-l/2,1-l/2,\dots,l/2\}
\end{align*} 
where $R^2=\alpha^2+\frac{4\beta\gamma q^{-1}}{(q-q^{-1})^2}$.
In particular $M_l$ has eigenvalue $0$ if and only if $l$ is even or
\begin{align*}
  0=\rho_r\rho_{-r}=-(q^r-q^{-r})^2\left(\alpha^2+\beta\gamma q^{-1}
                    \left(\frac{q^r+q^{-r}}{q-q^{-1}}\right)^2\right).
\end{align*}  
The second case is equivalent to $c=c(n)$ for some $n\in I_l$.
As this case is excluded by assumption $q^{-l}\in J^c$ if and only if $l$ is
even.

Let $\C w$ denote the one-dimensional representation of $\UslZ$ uniquely
determined by
$E{\cdot}w=0$, $F{\cdot}w=0$, $K{\cdot}w=-w$.
By means of a base change the matrix (\ref{matrix}) corresponding to $-q^{-l}$
can be transformed into the matrix of the transpose of the left $X_c$-action
on the finite
dimensional $\UslZ$-module $\C w\otimes V_l$. The eigenvalues of this action
can be computed by means of \cite[Prop. 4.6]{a-MullSch}. In particular
the $X_c$-action has a nontrivial kernel if and only if
\begin{align*}
  0=(\rho_r+2\alpha)(\rho_{-r}+2\alpha)=(q^r+q^{-r})^2
     \left(\alpha^2-\beta\gamma q^{-1}
                    \left(\frac{q^r-q^{-r}}{q-q^{-1}}\right)^2\right)
\end{align*}  
for some $r\in I_l$. This equation is equivalent to
$c=\frac{1}{(q^r-q^{-r})^2}$, $r\neq 0$ or $c=\infty$, $r=0$. Notice that
$(q^r-q^{-r})^2-c(n)^{-1}=(q^{r+n}{+}
q^{-(r+n)})(q^{r-n}{+}q^{-(r-n)})\neq 0$ for all
$r,n\in \N_0/2$ and therefore these cases are not excluded.
\end{proof}

\section{Differential Calculus over $\podl$}
For the convenience of the reader the notion of differential calculus from
\cite{a-Woro2} is recalled. 
A \textit{first order differential calculus} (FODC)
over an algebra $\B$ is a $\B$-bimodule $\Gamma$ together with a
$\C$-linear map
\begin{equation*}
  \dif:\B\rightarrow\Gamma
\end{equation*}
such that $\Gamma=\Lin_\C\{a\,\dif b\,c\,|\,a,b,c\in\B\}$ and $\dif$
satisfies the Leibniz rule
\begin{align*}
  \dif(ab)&=a\,\dif b + \dif a\,b.
\end{align*}    
Let in addition $\cA$ denote a Hopf algebra and
$\Delta_\B:\B\rightarrow \B\otimes\cA$ a right $\cA$-comodule algebra
structure on $\B$.
If $\Gamma$ possesses the structure of a right $\cA$-comodule
\begin{equation*}
  \Delta_\Gamma:\Gamma\rightarrow\Gamma\ot \cA
\end{equation*}
such that
\begin{equation*}
\Delta_\Gamma(a\dif b\,c)=(\Delta_\B a)((\dif\otimes\id)\Delta_\B b)
                          (\Delta_\B c)
\end{equation*}
then $\Gamma$ is called \textit{right covariant}.
A FODC $\dif:\B\rightarrow\Gamma$
over $\B$ is called \textit{inner} if there exists an element
$\omega\in\Gamma$ such that $\dif x=\omega x-x\omega$ for all $x\in\B$.
For further details on first order differential calculi consult
\cite{b-KS}.

Let $U $ denote a Hopf algebra with bijective antipode and $L \subset U $
a left coideal subalgebra, i.e.
$\kow_L :L \rightarrow \U \ot L $. Consider a
tensor category $\mathcal{C}$ of finite dimensional left $U $-modules. Let
$\A :=U ^0_{\mathcal{C}}$ denote the dual Hopf algebra generated by the
matrix coefficients of all $U $-modules in $\mathcal{C}$.
Assume that $\A$ separates the elements of $\U$.
Define a right coideal subalgebra $\B \subset \A $ by
\begin{align}\label{Bdef}
  \B :=\{b\in \A \,|\,\pair{u}{b_{(1)}}b_{(2)}=0 \quad
  \text{for all $u\in L^+$}\},
\end{align}
where $L^+=\{u-\vep(u)\,|\,u\in L\}$.
Assume $L$ to be $\mathcal{C}$-semisimple, i.e. the restriction of any
$\U$-module in $\mathcal{C}$ to the subalgebra $L\subset \U$ is isomorphic
to the direct sum of irreducible $L$-modules. By \cite{a-MullSch}
Theorem 2.2 (2) this implies that $\A $ is a faithfully flat $\B $ module.

In this situation right covariant first order differential calculi
over $\B$ can be classified via certain left ideals of
$\B^+$ \cite{a-Herm01}. More explicitly the subspace
\begin{align}
  \lid=\Big\{\sum_i a_i^+\vep(b_i)\,\Big|\,\sum_i\dif a_i\,b_i=0\Big\}
  \subset \B^+
\end{align}  
is a left ideal which determines the differential calculus uniquely.
To this left ideal one associates the vector space
\begin{align*}
T^\vep=\{f\in \B^\circ\,|\, f(x)=0 \textrm{ for all }x\in\lid\}
\end{align*}
and the so called \textit{quantum tangent space}
\begin{align*}
  T=(T^\vep)^+=\{f\in T^\vep\,|\, f(1)=0\}.
\end{align*}  
The dimension of a first order differential calculus is defined by
\begin{align*}
  \dim \Gamma=\dim_\C \Gamma/\Gamma\B^+=\dim_\C\B^+/\lid.
\end{align*}  

\begin{proposition}\label{corresp}{\em \cite[Cor.~1.2]{a-HeKo01p}}
There is a one to one correspondence between
$n$-dimensional covariant FODC over $\B $ and $(n+1)$-dimensional
subspaces $T^\vep\subset \B ^\circ $ such that 
\begin{align}
\vep \in T^\vep,\quad \Delta T^\vep \subset \B ^\circ \ot T^\vep,\quad T^\vep
 L\subset T^\vep.
\end{align}
\end{proposition}
A covariant FODC $\Gamma$ over $\B$ is called
\textit{irreducible} if it does not possess any nontrivial quotient
(by a right covariant $\B$-bimodule). Note that this property is equivalent
to the property that $T^\vep_\Gamma$ does not possess any 
right $L$-invariant left $\B^\circ$-subcomodule $\tilde{T}$ such that
$\C\cdot\vep\varsubsetneq\tilde{T} \varsubsetneq T^\vep_\Gamma$.

For a family of right covariant FODC
$(\Gamma_i,\dif_i)_{i=1,\dots,k}$ define $\dif=\oplus_i \dif_i:\B \rightarrow
\oplus_i\Gamma_i$. Then
$\Gamma=\B\dif \B\subset \oplus_i \Gamma_i$ is a covariant FODC with
differential $\dif$
which is called the \textit{sum} of the calculi $\Gamma_1,\dots,\Gamma_k$
\cite{a-HeckSchm2}. The left ideal corresponding to $\Gamma$ is given by
$\lid_\Gamma=\cap_i \lid_{\Gamma_i}$ and therefore the relation
$T_\Gamma=T_{\Gamma_1}+\dots+T_{\Gamma_k}$ of quantum tangent spaces holds.
A sum of covariant differential calculi is called a \textit{direct sum} if
$\Gamma=\oplus_i \Gamma_i$ is a direct sum of bimodules.
This condition is equivalent to $T_\Gamma=\oplus_i T_{\Gamma_i}$.

As an immediate consequence of Corollary \ref{corresp}, 
Lemma \ref{equiv} and Proposition \ref{locfin} one obtains the following
classification result for differential calculi over $\podl$.  
\begin{theorem}\label{podlclasstheo}
Assume $0\neq c\in J_2$. For $\lambda\in J^c$ let $\Gamma_{\lambda}$ denote
the uniquely determined covariant FODC over $\podl$ such that
$T_{\Gamma_\lambda}^\vep=V_\lambda+ \C\vep$.
Then $\Gamma_\lambda$ is irreducible and any finite
dimensional covariant FODC $\Gamma$ over $\podl$ is isomorphic to a direct
sum
\begin{align*}
  \Gamma=\bigoplus_{\lambda\in J} \Gamma_\lambda
\end{align*}  
for some finite subset $J\subset J^c$.
\end{theorem}  

For any covariant FODC $\Gamma$ with corresponding left ideal $\lid$ and
quantum tangent space $T$ consider the projection
\begin{align*}
  P_r:\Gamma\ot_\B\A\rightarrow \Gamma\ot_\B\A,\quad
  \gamma\ot a\mapsto \gamma_{(1)}\ot S(\gamma_{(2)})\vep(a)
\end{align*}    
onto the subspace $(\Gamma\ot_\B\A)_{inv}\subset \Gamma\ot_\B\A$ of right
coinvariant elements. The relation $\dif b\ot a=\dif(b_{(1)})\ot
S(b_{(2)})b_{(3)}a$ implies that the right $\A$ module $\Gamma\ot_\B\A$
is generated by the elements $P_r(\dif b\ot 1)$, $b\in\B$.
For any $a=\sum_i a_i^+\vep(b_i)\in\lid$ where $\sum_i\dif a_i b_i=0$
one obtains
\begin{align*}
    P_r(\dif a\ot 1)=P_r\left(\sum_i \dif a_i\ot b_i\right)=0.
\end{align*}
Therefore $P_r$ induces a well defined surjection
\begin{align}\label{surj}
  \B^+/\lid\rightarrow (\Gamma\ot_\B\A)_{inv},\quad
  b\mapsto P_r(\dif b\ot 1).
\end{align}

\begin{lemma}\label{pairing}
  The pairing
\begin{align*}
  (\Gamma\ot_\B\A)_{inv}\times T\rightarrow \C,\quad
  (\dif a\ot b,X)\mapsto X(a)\vep(b)
\end{align*}
is non-degenerate. Further $b\in\lid$ if and only if $b\in \B^+$ and
$P_r(\dif b\ot 1)=0$.
\end{lemma}
\begin{proof}
  To verify the first statement note that by construction the elements
  $P_r(\dif b\ot 1)$, $b\in\B$ separate $T$. On the other hand
  (\ref{surj}) implies
  $\dim_\C((\Gamma\ot_\B\A)_{inv})\le\dim_\C\B^+/\lid=\dim_\C T$ and
  therefore $T$ separates $(\Gamma\ot_\B\A)_{inv}$ and (\ref{surj}) is an
  isomorphism.
\end{proof}  

\begin{lemma}\label{freiekalk}
  Let $W\subset \B$ be a right $\A$-subcomodule then $\dif W$ generates
  $\Gamma$ as a right $\B$-module if and only if the elements of $W$
  separate the quantum tangent space $T_\Gamma$.
  If $\dim W=\dim \Gamma$ and the elements of $W$ separate $T_\Gamma$
  then $\Gamma$ is a free right $\B$-module generated by the
  differentials of an arbitrary basis of $W$.
\end{lemma}  

\begin{proof}
  Let $\Gamma'\subset \Gamma$ denote the right $\B$-module generated by
  $\dif W$. Then as $\A$ is a faithfully flat left $\B$-module
  \begin{align*}
    \Gamma'=\Gamma\Longleftrightarrow \Gamma'\ot_\B\A=\Gamma\ot_\B\A
    \Longleftrightarrow (\Gamma'\ot_\B\A)_{inv}=(\Gamma\ot_\B\A)_{inv}.
  \end{align*}
  Now, if $W$ separates $T_\Gamma$ then $(\Gamma'\ot_\B\A)_{inv}$
  separates $T_\Gamma$ and therefore by Lemma \ref{pairing}
  coincides with $(\Gamma\ot_\B\A)_{inv}$.
  Conversely, if $\Gamma'=\Gamma$ then $(\Gamma'\ot_\B\A)_{inv}$ separates
  $T_\Gamma$ and therefore the elements of $W$ separate $T_\Gamma$.
  This proves the first statement.

  To prove the second statement let $\Gamma^{''}$ denote the
  free right $\B$-module generated by the differentials of an arbitrary basis
  $e_1,\dots,e_k$ of $W$. Then, as above,
  \begin{align*}
    \Gamma^{''}\cong\Gamma
    &\Longleftrightarrow \Gamma^{''}\ot_\B\A\cong\Gamma\ot_\B\A\\
    &\Longleftrightarrow P_r(\dif e_i\ot 1), i=1,\dots,k
                         \textrm{ form a basis of }(\Gamma\ot_\B\A)_{inv}.
  \end{align*}
  In view of Lemma \ref{pairing} this poperty is equivalent to the
  nondegeneracy of the pairing between $W$ and $T_\Gamma$.
\end{proof}

Combining the above Lemma with Theorem \ref{podlclasstheo} one can now
classify all covariant FODC over $\podl$ generated as right $\podl$-modules
by the differentials $\dif e_i$, $i=-1,0,1$. The straightforward 
calculations of the pairing of the tangent spaces with the generators
$e_i$, $i=-1,0,1$, are omitted.
\begin{cor}
  For $c\in J_2\setminus\{0,\infty, (q^{1/2}-q^{-1/2})^{-2}\}$
  there exists exactly one covariant FODC $\Gamma_{q^{-2}}$ over
  $\podl$ which is generated by
  $\{\dif e_i\,|\,i=-1,0,1\}$ as a right $\podl$-module. The elements
  $\{\dif e_i\,|\,i=-1,0,1\}$ form a right $\podl$-module basis of this
  calculus.

  For $c=\infty$ there exist exactly three covariant FODC over $\podl$ which
  are generated by $\{\dif e_i\,|\,i=-1,0,1\}$ as a right $\podl$-module.
  One of them, $\Gamma_{-1}$, is one-dimensional, the elements
  $\{\dif e_i\,|\,i=-1,0,1\}$ form a right $\podl$-module basis of each of
  the other two calculi $\Gamma_{\pm q^{-2}}$.

  For $c=(q^{1/2}-q^{-1/2})^{-2}$ there exist exactly two covariant FODC over
  $\podl$ which are generated by $\{\dif e_i\,|\,i=-1,0,1\}$ as a right
  $\podl$-module. One of them, $\Gamma_{-q^{-1}}$, is two-dimensional,
  the elements $\{\dif e_i\,|\,i=-1,0,1\}$ form a right $\podl$-module basis
  of the other calculus $\Gamma_{q^{-2}}$.
\end{cor}  
For generic value of $c$ the above corollary reproduces the results
obtained in \cite{a-ApelSchm94} by means of computer calculations.

The odd dimensional covariant FODC $\Gamma_{q^{-l}}$, $l\in2\N$, for arbitrary
$c$ and $\Gamma_{-q^{-l}}$, $l\in2\N$, for $c=\infty$ can be
explicitly constructed by a method by U.~Hermisson. To match the above
conventions the relevant lemma from \cite{a-Herm98}, \cite{a-Herm01} is
cited in terms of right comodule algebras.

Let $\A$ denote a coquasitriangular Hopf algebra with universal
$r$-form $\rform $ and $\B$ a right $\A$-comodule algebra. Let $\nu$ be a
comodule algebra endomorphism of $\B$. Let further $W\subset \B$ denote a
finite dimensional right $\A$ subcomodule and $W'=\Hom(W,\C)$ the dual
comodule defined by $(\kow f)(w)=(f\ot S^{-1})\kow w$ for $w\in W$, $f\in W'$.
More explicitly, if $\{b_1,\dots,b_N\}$ is a basis of $W$ and
$\{\gamma^1,\dots,\gamma^N\}\subset W'$ the dual basis then
$\kow b_i=b_j\ot \psi^j_i$ implies $\kow \gamma^i=
\gamma^j\ot S^{-1}(\psi^i_j)$.

\begin{lemma}\label{ulrich}
The free right $\B$-module $W'\ot \B$ can be endowed with a right
$\A$-covariant $\B$-bimodule structure by
\begin{align*}
  afb:=f_{(0)}\nu (a_{(0)})b\,\rform (a_{(1)},f_{(1)}),\qquad
  a,b\in \B,f\in W'
\end{align*}
and will be denoted by $\Gamma _{\rform,\nu,W}$.
Moreover, if $\omega=\sum_{i=1}^N\gamma^i b_i\in \Gamma _{\rform,\nu,W} $
denotes the canonical invariant element then
$\dif:\B\to \Gamma _{\rform,\nu,W}$, $\dif b:=\omega b-b\omega$,
defines a covariant FODC $(\dif \B\cdot \B,\dif )$ over $\B$.
\end{lemma}

\begin{lemma}
  The quantum tangent space of the differential calculus $\Gamma$ described
  in Lemma \ref{ulrich} is the linear span of the functionals
  $\chi_i\in\B^\circ$, $i=1,\dots,N$, defined by
  \begin{align*}
     \chi_i(a)=\rform (\nu(a),S^{-1}(b_i))-\vep(b_i)\vep(a).
  \end{align*}  
\end{lemma}  

\begin{proof}
  For $a\in\B$ one obtains
  \begin{align*}
     -P_r(\dif a\ot 1)&=P_r\left(\sum_i(a\gamma^i b_i-\gamma^i b_i a)\ot 1
                        \right)\\
                     &=\sum_{i,j}\gamma^j\ot S^{-1}(\psi^i_j)\left[
                       \rform (\nu(a),S^{-1}(b_i))-\vep(b_i)\vep(a)\right].
  \end{align*}
  Therefore by Lemma \ref{pairing}
  \begin{align*}
    a\in\lid&\Longleftrightarrow a\in\B^+\mbox { and } P_r(\dif a\ot 1)=0\\
    &\Longleftrightarrow a\in\B^+\mbox{ and }\chi_i(a)=0 \quad
    \forall i=1,\dots, N.
  \end{align*}  
\end{proof}  

Let $V(n)$, $n\geq 1$, denote the $(2n{+}1)$-dimensional $\UslZ$-submodule of
$\podl$ with highest weight vector $b_1=e_1^n$.
For $\nu=\id$ the quantum tangent space $T$ of the differential calculus
$\Gamma$ from Lemma \ref{ulrich} satisfies
\begin{align*}
  \rform(\nu(\cdot),S^{-1}(e_1^n))=\chi_1(\cdot)+\vep(e_1^n)\vep(\cdot)\in
  T^\vep=T\oplus \C\vep.
\end{align*}
The standard universal $r$-form of $\SLZ$ is defined by
\begin{align*}
  \rform(u^i_j\otimes u^k_l)=q^{-1/2}\begin{cases}q&\textrm{ if }i=j=k=l\\
                                     1&\textrm{ if }i=j\neq k=l\\
                                     q-q^{-1}&\textrm{ if } j=k<i=l\\
                                     0&\textrm{else. } 
    \end{cases}
\end{align*}  
In particular $\rform(a\ot u^2_1)=0$ for all $a\in\SLZ$ and therefore
$\bar{\chi}:\B\to\C$,
\begin{align*}
\bar{\chi}(a):=\vep(e_1)^{-n}\rform(a,S^{-1}(e_1^n))
              =\vep(e_1)^{-n}\rform(S(a),e_1^n)
              =\rform(S(a),(u^1_1)^{2n})=\rform(a,(u^2_2)^{2n})
\end{align*}  
is a character which satisfies
\begin{align*}
  \bar{\chi}(e_i)=q^{-2ni}\vep(e_i).
\end{align*}
Since $\bar{\chi}=\psi^0_{q^{-2n}}\in V_{q^{-2n}}$
and $\dim V_{q^{-2n}}\oplus \C\vep=2n+2\ge \dim T^\vep $ one obtains
$T^\vep=V_{q^{-2n}}\oplus \C\vep$.
Hence the differential calculus $\Gamma$ coincides with
$\Gamma_{q^{-2n}}$. Similarly the differential calculus $\Gamma_{-q^{-l}}$,
$l\in 2\N$, over $\mathcal{O}_q(\mathbb{S}^2_\infty)$ can be realized using
the comodule algebra endomorphism $\nu:e_i\mapsto -e_i$.

Note that
\begin{align}\label{inv}
  (\Gamma_{q^{-2n}}\ot_\B\A)_{inv}=(\Gamma_{\rform,\ids,V(n)}\ot_\B\A)_{inv}
\end{align}  
where as above $V(n)$ denotes the $(2n{+}1)$-dimensional representation of
$\UslZ$. Indeed, by the above remarks $\Gamma=\Gamma_{q^{-2n}}$ can be
considered as a right $\B$-submodule of $\Gamma_{\rform,\ids,V(n)}$ and as
$\A$ is a flat $\B$-module this implies
$\Gamma_{q^{-2n}}\ot_\B\A\subset \Gamma_{\rform,\ids,V(n)}\ot_\B\A$.
As $\dim (\Gamma_{\rform,\ids,V(n)}\ot_\B\A)_{inv}=2n+1$ by construction and
$\dim(\Gamma_{q^{-2n}}\ot_\B\A)_{inv}=2n+1$ by Lemma \ref{pairing} the
identification (\ref{inv}) follows.
Now (\ref{inv}) implies
$\Gamma_{q^{-2n}}\ot_\B\A=\Gamma_{\rform,\ids,V(n)}\ot_\B\A$ and by faithful
flatness of $\A$ this in turn gives
$\Gamma_{q^{-2n}}=\Gamma_{\rform,\ids,V(n)}$.
Thus one has the following proposition.

\begin{proposition}\label{freiinner}
For any $0\neq c\in J_2$ the FODC $\Gamma_{q^{-2n}}$, $n\in \N$, is
isomorphic to $\Gamma _{\rform,\ids,V(n)}$. For $c=\infty$ the FODC 
$\Gamma_{-q^{-2n}}$ is isomorphic to $\Gamma _{\rform,\nu,V(n)}$
where $\nu(e_i)=-e_i$. In particular $\Gamma_{\pm q^{-2n}}$ are free left and
right $\B$ modules and inner first order differential calculi.
\end{proposition}

\begin{remark}
  Covariant FODC over $\podlo$ are
  qualitatively different from those over $\podl$, $0\neq c\notin J_2$.
  Let $\tilde{\Gamma }_{kl}$ denote the $kl+k+l$-dimensional FODC over
  $\podlo$ with quantum tangent space
  $\tilde{T}_{kl}=
  \Lin_\C\{E^iF^j\,|\,0\leq i\leq k,0\leq j\leq l,(i,j)\neq (0,0)\}$.
  By Proposition \ref{corresp} and \cite[Lemma 5.3]{a-HeKo01p}
  any covariant FODC over $\podlo$ can be written as a (not necessarily
  direct) sum of calculi $\tilde{\Gamma }_{kl}$ for certain $k,l$.
  In particular the only irreducible calculi are $\tilde{\Gamma }_{10}$
  and $\tilde{\Gamma }_{01}$ constructed in \cite{a-Kolb01p}.  
\end{remark}

\bibliographystyle{amsalpha}
\bibliography{litbank2}

\end{document}